\documentclass{article}
\usepackage{graphicx} % Required for inserting images
\usepackage[utf8]{inputenc}
\usepackage{color}
\usepackage{amssymb}

\usepackage{bm}  
\usepackage{hyperref}       % hyperlinks
\usepackage{fullpage}       % hyperlinks
\usepackage{url}            % simple URL typesetting
\usepackage{booktabs}       % professional-quality tables
\usepackage{amsfonts}       % blackboard math symbols
\usepackage{amssymb}        % blackboard math symbols
\usepackage{nicefrac}       % compact symbols for 1/2, etc.
\usepackage{microtype}      % microtypography
\usepackage{tikz}           % for graph figures
\usepackage{centernot}

\usepackage[textsize=tiny]{todonotes}
\usepackage[normalem]{ulem}
\newcommand{\stkout}[1]{\ifmmode\text{\sout{\ensuremath{#1}}}\else\sout{#1}\fi}

\usepackage{graphicx}
\graphicspath{{rawfigs/}}

\usepackage[ruled,vlined,linesnumbered]{algorithm2e}
\usepackage{amsmath}
\usepackage{amssymb}
\usepackage{amsthm}
\usepackage{bm}
\usepackage{graphicx}
\graphicspath{{figs/}}
\usepackage{mathtools}
\usepackage{subcaption}

\newcommand{\comment}[1]

\theoremstyle{remark}
\newtheorem*{remark}{Remark}

\newtheorem{theorem}{Theorem}[section]
\newtheorem{lemma}[theorem]{Lemma}
\newtheorem{corollary}{Corollary}
\newtheorem{defn}{Defintion}[section]
\newtheorem{conj}{Conjecture}[section]

\newtheorem{example}{Example}[section]

\numberwithin{equation}{section}

\renewcommand{\kappa}{\varkappa}
\newcommand{\be}{\begin{equation}}
\newcommand{\ee}{\end{equation}}

\newcommand{\bq}{\begin{eqnarray}}
\newcommand{\eq}{\end{eqnarray}}

\newcommand{\ba}{\begin{array}}
\newcommand{\ea}{\end{array}}

\title{Exact Reduction of the Generalized Lotka-Volterra Equations via Integral and
Algebraic Substitutions}
\author{Rebecca Morrison}

\date{January 2021}

\begin{document}

\maketitle

\abstract{Systems of interacting species, such as biological environments or chemical reactions, are often
    described mathematically by sets of coupled ordinary differential equations. While a large
    number $\beta$ of species may be involved in the coupled dynamics, often only $\alpha< \beta$ species are of
    interest or of consequence. In this paper, we explore how to construct models that include
    only those given $\alpha$ species, but still recreate the dynamics of the original $\beta$-species model.
    %This type of dimension reduction does not yield a model that is computationally easier to solve.
    Under some conditions detailed here, this reduction can be completed exactly, such that the
    information in the reduced model is exactly the same as the original one, but over fewer
    equations.  Moreover, this reduction process suggests a promising type of approximate model---no
    longer exact, but computationally quite simple.
\\

}
%Department of Computer Science, University of Colorado Boulder, Boulder, CO, USA
%\def\corrAuthor{Rebecca E. Morrison, 1111 Engineering Dr, Boulder, CO 80309}
%\def\corrEmail{rebeccam@colorado.edu}
%\extraAuth{}
%
\section{Introduction}\label{sec:int} Consider an environment in which a large number of species interact. This
could be, for example, a chemical reaction (with reacting chemical species) or an ecological
system (with interacting organisms). To model the behavior of interacting species in many such
environments, we often use the generalized Lotka-Volterra equations---a set of coupled
ordinary differential equations (ODEs) \cite{yasuhiro1996global}. The generalized Lotka-Volterra (GLV) equations extend
2-species predator-prey models common in undergraduate differential equations classes to any number
of species. They are in fact much more general than a predator-prey type system: they provide a
framework to describe the time-dynamics of any number $\beta$ of interacting species, allowing for
linear (growth rate) and quadratic (interaction) terms.

In many applications, a large number of species may play a role in the dynamical community behavior.
However, from a modeling perspective, the species of interest or of consequence may be restricted to
a much smaller subset of $\alpha$ species, where $\alpha < \beta$. This occurs in fields as diverse as ecology
\cite{wangersky1978lotka}, epidemiology \cite{dantas2018calibration,li1995global}, and chemical kinetics
\cite{williams2008detailed,jones1988global,frassoldati2009simplified}. 
%For example, the state-of-the-art model for methane combustion,
%GRI-Mech \cite{grimech}, includes over 50 interacting species. But, many simplified methane
%combustion models include only about four to ten species
%\cite{jones1988global,frassoldati2009simplified}. As another example, consider that a standard
%epidemiological model will include humans and the disease carrier, such as mosquitoes. Often these
%models differentiate ``sub-species'' of each, that is, populations of susceptible, exposed,
%infected, and recovered humans and carriers \cite{li1995global}. Clearly, these models omit many
%other ecological species (horses, cattle, etc.) and sub-populations (quarantined or hospitalized
%humans) that could be involved in this system. 
There are good reasons to omit many species. First, a modeler may not have access to data about
certain species' concentrations, which would be needed to define initial conditions or calibrate
additional model parameters. Second, the modeler may not actually know what species should be included.
Third, the more species included in the model, then in general the more computationally expensive
the model becomes. Except in very special (and usually unrealistic) cases, a system of GLV equations
will not admit a closed-form solution; instead we solve these systems with computational models.
Therefore, it is common to build reduced models which include only $\alpha<\beta$ given species.

One immediate question that arises in this context is the following: Given a system of $\beta$ species,
suppose only $\alpha$ are known, or of interest. What is the best reduced deterministic model, in terms
of only those given $\alpha$ species? Of course, to answer this question we must define what is meant by
``best.'' To do so, let us first consider the landscape of computational models---one of science,
mathematics, and computation. We will broadly classify this consideration into three major areas:
model validation, model reduction, and computational implementation, and what the major questions
are in each.
\begin{itemize}
    \item[I] \emph{Model validation.} Does the mathematical model adequately represent the
        scientific system in question?
    \item[II] \emph{Model reduction.} How much error is incurred by use of the reduced model
        compared to the original, or high-fidelity, model?
    \item[III] \emph{Computational implementation.} Is the reduced model less computationally expensive
        than the original model, and by how much?
\end{itemize}
(Point III is closely related to model verification, a process which checks that the computational
model solves the mathematical model correctly \cite{oberkampf2010verification}.)
%In this paper, all computational implementations are %assumed correct and verified.)
Returning to the question of the best reduced model: the best reduced model would address all three areas,
\emph{i.e.,} well-represent the scientific system
under study (I); recreate the dynamics of the high-fidelity model (II); and be computationally easier
to solve than the high-fidelity model (III).

The first topic, model validation, is a recently growing and still very open field. For a general
description, see, for example, \cite{oliver2015validating}. For a few specific works, see
\cite{bayarri2007framework,morrison2018representing}. In this paper, the original model of $\beta$
species represents the true (physical, chemical, biological, etc.) system under study.

In many types of model reduction of systems of differential equations, the goal is to reduce the
computational cost (III) while controlling the incurred error (II). This certainly makes sense as an
objective---one would expect a reduced model to offer computational savings over the original model.
In order to gain computationally efficiency, these types of reductions are not exact but instead
offer an approximation of the species behavior. For example, eigendecompositions may yield an
approximation of the static equilibrium state and not the dynamical behavior.  Other techniques
reduce the computational complexity while maintaining time-dynamics, such as volume averaging
\cite{murdoch1994continuum}, perturbation theory \cite{gear2003projective}, spectral analysis
\cite{mezic2005spectral} and separation of fine and coarse scale variables
\cite{tartakovsky2011dimension}. Related work identifies the exact dynamics on the relevant slow
manifold, using, for example, fast-slow decompositions \cite{haller2017exact}, or polynomial
approximations up to a specified accuracy \cite{kazantzis2010new}.  Identification of symmetries
can reveal reduced exact dynamics in initial value problems \cite{zhdanov2002higher} and the
KdV--Zakharov--Kuznetsov equation \cite{sahoo2017lie}. For a good overview of more methods, see
\cite{pavliotis2008multiscale}.

%However, these common model reduction techniques may be undesirable or unused. A major reason for
%this is that some of these techniques will output quantities that do not directly correspond to
%individual species. For example, PCA will \dots
%Some existing model reduction techniques and the ones presented in this work differ in what the term
%``reduction'' actually means. In many of the techniques above, model reduction implies a resulting
%model which is significantly easier to compute. 

In contrast, here our goal is to reduce the number of coupled equations that make up the model, or
equivalently, the number of species involved in the model. Variables  are eliminated in a way that
is exact, that is, \emph{without loss of information}. In terms of the points above, (I) the model
is an exact representation of the true system, and (II) the reduction incurs zero error. We
investigate two possibilities for this type of dimension reduction in the context of the generalized
LV equations, called \textbf{integral} and \textbf{algebraic} substitution. There are two defining
characteristics of these methods. First, they preserve the correspondence between the set of $\alpha$
species of interest as they appear in the original model and the resulting set after the reduction
occurs.  This property is termed \emph{species correspondence}, or simply \emph{correspondence}.
Second, they create a reduced model which contains the exact same information as the original
one, but with fewer equations. Although the resultant model is not necessarily better in a
computational sense, this method reveals a path towards model reduction that preserves
correspondence, closely approximates the dynamics, and is computationally quite simple.

%In other words, we aim to perform model reduction exactly (II) but without concern, at first, for the
%computational implementation (III).  This may not seem a practical goal. However, the process
%provides insight into the coupled and complex nature of these GLV equations, and later suggests
%further techniques that do both. That is, a desirable method reduces both the number of
%equations and their computational complexity, while still preserving species
%correspondence and dynamics, thereby achieving all goals (I - III).

The two methods discussed here do not automatically find the $\alpha$
species comprising the reduced set, but rather assume them given. Some techniques do choose this set
as a step within the reduction process itself. Doing so may favor species with high relative
concentrations, or those with slow dynamics, for example. But consider a combustion model in which we
are concerned about trace amounts of a contaminant, or an ecological model built to track  a species
near extinction. The methods presented here allow the modeler to include these critical species \emph{a priori}.

The rest of the paper is structured as follows. Section~\ref{sec:glv} outlines the generalized
Lotka-Volterra equations. Section~\ref{sec:mr} presents our main results about two types of
exact dimension reduction. Based on these results, Section~\ref{sec:app} presents some possible
approximate methods, and we conclude with Section~\ref{sec:con}.

%***************************************************************************************************
\section{Background: The generalized Lotka-Volterra equations}\label{sec:glv} Let us briefly
%\section{Results}\label{sec:glv} Let us briefly
review the GLV equations. Let $\bm{x}$ be the $\beta$-vector of species concentrations. Here, units refer to
the number of specimens per unit area, but specific units are omitted for this paper. In the framework of
generalized Lotka-Volterra equations, the system of ODEs is given as \begin{equation}
\frac{d\bm{x}}{dt} = \text{diag}(\bm{x})(\bm{b} + \bm{A}\bm{x}), \end{equation} where the $\beta$-vector
$\bm{b}$ is the intrinsic growth rate vector, and $\bm{A}$ is the $\beta\times \beta$ interaction matrix. Note that
there is one differential equation for each species variable. In particular, if equilibrium is
achieved, then the equilibrium solution is given by
\begin{equation} \bm{x} = -A^{-1}\bm{b}. \end{equation}
Moreover, there is coexistence, or a feasible equilibrium, if, at equilibrium,
$x_i > 0\,\, \forall i \in 1,\dots,\beta$.

Here it is assumed all $x_i(t) >0, t \geq 0$. Conditions on the entries of $A$ and $\bm{b}$ that
guarantee coexistence are given in \cite{barabas2016effect}; the authors study the effects of intra-
(within a species) and inter- (among different species) specific competition on coexistence of large
systems.  In \cite{grilli2017feasibility}, the authors explore conditions for stability under
perturbations, asymptotic feasibility, and equilibrium of large systems.  

%A significant amount of theoretical ecology is understood about these mathematical systems; we
%present just a bit of vocabulary and the corresponding literature here. 

%***************************************************************************************************
\section{Exact dimension reduction}\label{sec:mr} 
We begin with some definitions. Let $\alpha, \beta \in \mathbb{Z}^+$.
\begin{defn}
[$(\beta,\alpha)$-reducible] A system of $\beta$ differential equations that can be
converted to a set of $\alpha$ differential equations, where $\alpha < \beta$ and without loss of
information, is called \emph{($\beta,\alpha$)-reducible}. \end{defn}

\begin{defn}[Integral substitutions] Integral substitutions (IS) may allow for reductions from a system of
    $\beta$ equations to $\alpha$ by eliminating the variables $x_{\alpha+1}, \dots x_\beta$. During
    this process, we introduce the entire history, or memory, of a subset of the variables $x_1, \dots
    x_\alpha$.  \end{defn}
This approach is similar to the Mori-Zwanzig method of model reduction \cite{givon2004extracting}. 

\begin{defn}[Algebraic substitutions] Algebraic substitutions (AS) may allow for reductions from a system of
    $\beta$ equations to $\alpha$ by eliminating the variables $x_{\alpha+1}, \dots x_\beta$. During
    this process, we introduce higher-order derivatives of a subset of the variables $x_1, \dots
    x_\alpha$. \end{defn}
This process is similar to the algebraic reduction presented in \cite{harrington2017reduction}.
%Recall that the term ``reduction'' is used to mean a reduction in the number of species, or
%equivalently coupled equations, in the model. 
%Let us first examine $(\beta,s)$-reducibility in the simplest case: $\beta=2$ and $\alpha=1$.

%***************************************************************************************************
\begin{example}{The GLV equations are $(2,1)$-reducible via IS.}  \end{example} 
    %Note that $x_1$ only has time history; there is no spatial variation modeled here. 
With $\beta=2$, the original model is:
\begin{subequations}
\begin{align}
    \dot x_1 &=  b_1  x_1 + ( a_{11} x_1 +  a_{12} x_2) x_1 \label{eq:dx1}\\
    \dot x_2 &=  b_2  x_2 + ( a_{21} x_1 +  a_{22} x_2) x_2 \label{eq:dx2}.
\end{align}
\end{subequations}
The goal is to rewrite $x_2$ in terms of $x_1$ in equation \ref{eq:dx1}, specifically in the term
$a_{12} x_2 x_1$.

First, rearranging equation \ref{eq:dx1} gives
    \begin{align}x_2 &= \frac{1}{a_{12}}\left(\frac{1}{x_1}(\dot x_1 - b_1 x_1)  -
    a_{11}x_1\right) \equiv \bm{y_2^1} \label{eq:y21}.\end{align}
    We denote this quantity $\bm{y_2^1}$ because it is a representation of $x_2$ that only depends on
    $x_1$. With regards to notation, bold type indicates any such introduced variables to more easily distinguish them from
    the $x_i$s. Also, a subscript indicates which variable the new one is replacing, and the
    superscript shows which variables this new one actually depends on. Now, substituting $\bm{y_2^1}$ back into \ref{eq:dx1}
     yields $0=0$.  Instead, let's substitute $\bm{y_2^1}$ for $x_2$ in \ref{eq:dx2}:
    \begin{equation}
    \dot x_2 =  b_2  \bm{y_2^1} + ( a_{21} x_1 +  a_{22} \bm{y_2^1}) \bm{y_2^1} .\label{eq:dy2}\end{equation}
    Integrating, we have,
\begin{align}
    x_2 &=  \int_0^t b_2  \bm{y_2^1} + ( a_{21} x_1 +  a_{22} \bm{y_2^1}) \bm{y_2^1}\equiv
    \bm{\chi_{2}^{1}} \label{eq:int2} .
\end{align}
Similarly, the symbol $\bm{\chi_2^1}$ means this is a variable that is equivalent to $x_2$ but only in terms
of $x_1$. Finally, \ref{eq:dx1} becomes:
\begin{equation}
    \dot x_1 =  b_1  x_1 + ( a_{11} x_1 +  a_{12} \bm{\chi_2^{1}} )x_1 .\label{eq:chi1}
\end{equation}
We now have a system of a single differential equation, in terms of $x_1$ and its memory.
%as represented by the integral equation \ref{eq:int2}. 
Note that this process has preserved species correspondence and that there is no loss of
information: the variable $x_1$ and its derivative in equation~\ref{eq:chi1} are equivalent to that
in \ref{eq:dx1}.

%***************************************************************************************************
\begin{example}{The GLV equations are also $(2,1)$-reducible via AS.}
\end{example}
The first step here is the same as that of the previous subsection, yielding $\bm{y_2^1}$.
Next, however, by differentiating $\bm{y_2^1}$, we have:
    \begin{equation}
        \bm{\dot{y}_{2}^{1}} = \frac{1}{a_{12}} \left( \frac{- \dot x_1}{x_1^2} \left( \dot x_1 - b_1
        x_1\right) +
        \frac{1}{x_1}\left(\ddot x_1 - b_1 \dot x_1\right) - a_{11}\dot x_1 \right). \label{eq:dy21}
    \end{equation}
Equations \ref{eq:y21} and \ref{eq:dy21}  express $x_2$ and $\dot x_2$, respectively, in terms of $x_1$, $\dot x_1$ and $\ddot
x_1$. Finally, we can rewrite the second equation as
\begin{equation} \bm{\dot y_2^1} = b_2 \bm{y_2^1} + (a_{21}x_1 + a_{22}\bm{y_2^1})
\bm{y_2^1}.\label{eq:alg1} \end{equation}
Again, equation~\ref{eq:alg1} is a system of a single differential equation, but here in terms of $x_1$ and its
derivatives, $\dot x_1$ and $\ddot x_1$.

\begin{remark}[Resulting functional form]\label{rem:fun}
Both the integral and algebraic forms, manipulated in this way, yield a single
differential equation in terms of $x_1$. Both methods respect species correspondence and are exact.
A major difference between the two methods is revealed by inspection of the structure of the
resulting sole equations. After reduction via the integral method, the structure of the final equation
resembles that of the initial equation \ref{eq:dx1} for $x_1$. The functional form matches in the
placement of the variables $x_1$ and $\bm{\chi_2^{1}}$ (in place of $x_2$), and also the remaining
constants $b_1$, $a_{11}$, and $a_{12}$ (which do not appear in the second equation of the original
system).  In contrast, after reduction via the algebraic method, the resultant equation has the
structure of \ref{eq:dx2}. This suggests that in an applied setting,  one type of reduction may be
advantageous, depending on what information is known about the high-fidelity model, such as various model
parameters.\end{remark}

For the remainder of the paper, we will focus on IS, as opposed to AS. The techniques and results
between the two methods are quite similar. An advantage
of IS is its notational convenience (as described in Remark~\ref{rem:fun}): the
reduced equations preserve the structure of the original ones, which simplifies tracking the
reductions, in theory.
%That is, in the example above, the
%integral method results in a differential equation for $x_1$ that has the same structure as the
%original ODE for $x_1$. On the other hand, AS results in a differential equation
%that only involves variable $x_1$, but is more similar to the original equation for $x_2$.

\begin{remark}[$(\beta,\beta-1)-reducibility$]
In a similar way, we can reduce any generalized LV system of $\beta$ species to one of $\beta-1$ species.
Note that this level of reduction, from $\beta$ to $\beta-1$ equations, can happen when the
ODEs describe the dynamics of fractional concentrations. In that case, the extra constraint that
$\sum_{i=1}^\beta x_i = 1$ readily allows for a reduction to $\beta-1$ equations, since, for example, $x_\beta$
can be expressed as $x_\beta = 1 - \sum_{i=1}^{\beta-1} x_i$. In this work, however, there is no such
restriction on the concentrations, and yet such a reduction is always possible.
\end{remark}

Reductions from two to one equation and from $\beta$ to $\beta-1$ are similar.
\begin{lemma}[]{The GLV equations are ($\beta, \beta-1$)-reducible via IS.}\label{lem:memSS-1}\end{lemma}
    \begin{proof}
    The original model for $\beta$ species is:
\begin{subequations}
    \begin{align}
        \dot x_1 &= b_1 x_1 + (a_{11} x_1 + a_{12}x_2 + \dots + a_{1S}x_\beta) x_1 \label{eq:1}\\
        \dot x_2 &= b_2 x_2 + (a_{21} x_1 + a_{22}x_2 + \dots + a_{2\beta}x_\beta) x_2 \label{eq:2}\\
        &\vdots \nonumber \\
        \dot x_\beta &= b_\beta x_\beta + (a_{\beta1} x_1 + a_{\beta2}x_2 + \dots + a_{\beta\beta}x_\beta) x_\beta \label{eq:S}.
    \end{align}
\end{subequations}
At this point, we want to rewrite $x_\beta$ in terms of the remaining $\beta-1$ variables, but we could use
any of the first $\beta-1$ equations to do so. Without loss of generality, we choose the second to
last one:
    \begin{align}
        x_\beta &= \frac{1}{a_{\beta-1,\beta}} \Biggl(\frac{1}{x_{\beta-1}} (\dot x_{\beta-1} - b_{\beta-1}
        x_{\beta-1})  -  a_{\beta-1,1} x_1 -  a_{\beta-1,2}x_2 - \dots - a_{\beta-1,\beta-1}x_{\beta-1}\Biggr) \equiv \bm{y_\beta^{1:\beta-1}}.
    \end{align}
The colon notation in $\bm{y_\beta^{1:\beta-1}}$ signifies that this new variable is written in terms of variables $x_1$ through $x_{\beta-1}$. 
    
    Now, substituting $\bm{y_\beta^{1:\beta-1}}$ into \ref{eq:S} yields
    \begin{align}x_\beta &= \int_0^t \left(b_\beta \bm{y_\beta^{1:\beta-1}} + \left(a_{\beta1}x_1 + a_{\beta2}x_2 +
    \dots + a_{\beta\beta} \bm{y_\beta^{1:\beta-1}}\right) \bm{y_\beta^{1:\beta-1}}\right)\nonumber \equiv \bm{\chi_{\beta}^{1:\beta-1}}.\end{align}
    And finally, substituting $\bm{\chi_{\beta}^{1:\beta-1}}$ back into the first $\beta-1$ equations:
\begin{subequations}
    \begin{align}
        \dot x_1 &= b_1 x_1 + \left(a_{11} x_1 + a_{12}x_2 +
        \dots+a_{1\beta}\bm{\chi_{\beta}^{1:\beta-1}}\right) x_1 \label{eq:chi1S}\\
        \dot x_2 &= b_2 x_2 + \left(a_{21} x_1 + a_{22}x_2 +\dots+
        a_{2\beta}\bm{\chi_{\beta}^{1:\beta-1}}\right) x_2 \label{eq:chi2S}\\
        &\vdots \nonumber\\
        \dot x_{\beta-1} &= b_{\beta-1} x_{\beta-1} + \left(a_{\beta-1,1} x_1 + a_{\beta-1,2}x_2 +
        \dots+a_{\beta-1,\beta}\bm{\chi_{\beta}^{1:\beta-1}}\right) x_{\beta-1} \label{eq:chiS-1S}.
    \end{align}
\end{subequations}
Thus the set of $\beta$ ODEs is reduced to $\beta-1$ without loss of
information.
    \end{proof}

%***************************************************************************************************
%\subsection{Model reduction, $s < \beta- 1$}
%Ideally, the above process would extend to larger systems, \emph{i.e.,} reduce without loss of
%information a GLV system of $\beta$ equations to any $s$, where $s <\beta$. 
In general, a system of GLV equations is not ($\beta, \alpha$)-reducible because of the coupled-ness, or
``entanglement'' that occurs between the species in the reduced
set via the introduced $\bm{y}$ and/or $\bm{\chi}$ variables. To see this, consider the case of
$\beta=3$, $\alpha=1$. With both the integral and algebraic methods, the model cannot be exactly
reduced, unless one of the coefficients is zero, thereby breaking the entanglement.

\begin{lemma}[]{The GLV equations are $(3,1)$-reducible if only one $a_{ij} = 0, i \neq
j$.}\end{lemma} \label{lem:31}
    \begin{proof}
We begin the process to reduce the system of three equations to one by eliminating
the variables $x_2$ and $x_3$.  When $\beta=3$, the original model is:
\begin{subequations}
    \begin{align}
        \dot x_1 &= b_1 x_1 + (a_{11} x_1 + a_{12}x_2 + a_{13}x_3) x_1 \label{eq:31}\\
        \dot x_2 &= b_2 x_2 + (a_{21} x_1 + a_{22}x_2 + a_{23}x_3) x_2 \label{eq:32}\\
        \dot x_3 &= b_3 x_3 + (a_{31} x_1 + a_{32}x_2 + a_{33}x_3) x_3 \label{eq:33}.
    \end{align}
\end{subequations}
Repeating the process described in the previous section, we can rewrite either equation \ref{eq:31} or \ref{eq:32} for
$x_3$. Using equation \ref{eq:32}, then
    \begin{align}x_3 &= \frac{1}{a_{23}}\left(\frac{1}{x_2}(\dot x_2 - b_2 x_2)  -
    a_{21}x_1 - a_{22}x_2 \right) \equiv \bm{y_3^{1:2}}.\label{eq:y312}\end{align}
    Next, substituting $\bm{y_3^{1:2}}$ into \ref{eq:33},
    \begin{align}
        \dot x_3 &=  b_3  \bm{y_3^{1:2}} + ( a_{31} x_1 +  a_{32} x_2 + a_{33}\bm{y_3^{1:2}}) \bm{y_3^{1:2}}
        \label{eq:dy3}.
    \end{align}
    Integrating,
\begin{align}
    x_3 &=  \int_0^t b_3  \bm{y_3^{1:2}} + ( a_{31} x_1 +  a_{32} x_2 + a_{33}\bm{y_3^{1:2}}) \bm{y_3^{1:2}}
    \equiv \bm{\chi_{3}^{1:2}}\label{eq:int3}.
\end{align}
Now we can substitute $\bm{\chi_{3}^{1:2}}$ into equations \ref{eq:31} and
\ref{eq:32}:
\begin{subequations}
    \begin{align}
        \dot x_1 &= b_1 x_1 + (a_{11} x_1 + a_{12}x_2 + a_{13}\bm{\chi_{3}^{1:2}}) x_1 \label{eq:321}\\
        \dot x_2 &= b_2 x_2 + (a_{21} x_1 + a_{22}x_2 + a_{23}\bm{\chi_{3}^{1:2}}) x_2 \label{eq:322}.
    \end{align}
\end{subequations}
At this point the model has been reduced from three equations to two. Consider if we now tried to
remove another variable, say $x_2$. The first step would be to rewrite equation  \ref{eq:321} so
that $x_2$ is alone on the left-hand side (LHS). However, now that $\bm{\chi_{3}^{1:2}}$ has been
introduced, we cannot cleanly separate $x_2$ out of the right-hand side (RHS) since it is embedded
inside the integral term.

A similar problem occurs with the algebraic approach.

Now assume of the six off-diagonal interaction terms is set to
zero. Without loss of generality, set $a_{13} = 0$. Then equation \ref{eq:321} becomes
\begin{equation}
    \dot x_1 = b_1 x_1 + (a_{11} x_1 + a_{12}x_2 ) x_1 \label{eq:dx1r}
\end{equation}
and so,
\begin{align}x_2 &= \frac{1}{a_{12}}\left(\frac{1}{x_1}(\dot x_1 - b_1 x_1)  -a_{11}x_1\right)
\equiv \bm{y_2^{1}}. \end{align}
Substituting $\bm{y_2^{1}}$ into equation \ref{eq:322},
\begin{align}
    \dot x_2 &= b_2 \bm{y_2^{1}} + (a_{21} x_1 + a_{22} \bm{y_2^{1}} + a_{23}\bm{\chi_{3}^{1:2}})
    \bm{y_2^{1}}\\
     &= b_2 \bm{y_2^{1}} + (a_{21} x_1 + a_{22}\bm{y_2^{1}} + a_{23}\bm{\chi_{3}^{1}}) \bm{y_2^{1}}
     \label{eq:y21-again}.
\end{align}
The second line above appears after replacing $\bm{\chi_{3}^{1:2}}$ with $\bm{\chi_{3}^{1}}$. This
variable $\bm{\chi_{3}^{1}}$ is found by replacing any explicit dependence on $x_2$ with
$\bm{y_2^{1}}$.
Specifically,
\begin{align} \bm{\chi_{3}^{1}} &= \int_0^t b_3  \bm{y_3^{1:2}} + ( a_{31} x_1 +  a_{32} x_2 +
    a_{33}\bm{y_3^{1:2}})\bm{y_3^{1:2}} \\
    &= \int_0^t b_3  \bm{y_3^{1}} + ( a_{31} x_1 +  a_{32} \bm{y_2^{1}} + a_{33} \bm{y_3^{1}}) \bm{y_3^{1}},
\end{align}
where, to in turn find $\bm{y_3^{1}}$, we modify equation~\ref{eq:y312}:
\begin{align}
    \bm{y_3^{1:2}} &= \frac{1}{a_{23}}\left(\frac{1}{x_2}(\dot x_2 - b_2 x_2)  - a_{21}x_1 - a_{22}x_2
    \right)\\
    &= \frac{1}{a_{23}}\left(\frac{1}{\bm{y_2^1}} (\dot{\bm{y^1_2}} - b_2 \bm{y^1_2})  - a_{21}x_1 - a_{22}\bm{y_2^{1}}
    \right) \equiv \bm{y_3^{1}}.
\end{align}
Integrating line~\ref{eq:y21-again},
\begin{align}
    x_2 &= \int_0^t b_2 \bm{y_2^{1}} + (a_{21} x_1 + a_{22}\bm{y_2^{1}} + a_{23}\bm{\chi_{3}^{1}})
    \bm{y_2^{1}}
    \equiv \bm{\chi_{2}^{1}},
\end{align}
and finally, 
\begin{equation}
    \dot x_1 = b_1 x_1 + (a_{11} x_1 + a_{12}\bm{\chi_{2}^{1}} ) x_1. \label{eq:x1a13}
\end{equation}
\end{proof}

So, we can exactly reduce the system of three ODEs to one differential equation if we assume $a_{13}
= 0$. However, note that a nested set of integrals comprises the term $\bm{\chi_{3}^{1}}$.  Indeed,
such a system of one differential equation with nested integral terms as equation~\ref{eq:x1a13} would
likely be more difficult to solve than the original set of three coupled ODEs
shown in \ref{eq:31} - \ref{eq:33}.

\begin{remark}[Entanglement] Multiple
variables can be eliminated if enough information (memory terms, higher derivatives) is known about
the remaining variables. Lemma~\ref{lem:31} shows the necessity of the assumption that some interaction term be zero
in order to exactly reduce (in this demonstration, that $a_{13} = 0$). This
assumption begins to reveal
the limitations of such reductions via substitutions---in particular, how, after a single
substitution, we cannot escape the entanglement of the remaining species. \end{remark}

    At the same time, the
methods here give insight into how one might approximate the role of eliminated variables in
terms of the reduced set. For example, in the case of $\beta=2, \alpha=1$, consider an approximation
$\bm{\hat \chi_{2}^{1}}$ in terms of $x_1$ and this extra information: \begin{equation}
\bm{\chi_{2}^{1}} \approx \bm{\hat \chi_{2}^{1}}\left(x_1, \dot x_1, \ddot x_1,
K(x_1),\dots\right),\end{equation} where $K(x_1)$ represents some memory kernel of $x_1$. We explore
this type of approximation in Section~\ref{sec:app}.
%We return to this topic in Section~\ref{sec:app}.
%(which is even more strikingly evident here than after just one variable removal). 

%With these reasons in mind, 
But first, we complete the reduction via IS for general $\beta$ and $\alpha$. Of course, in this
setting, a larger collection of interaction terms must be assumed zero.  Define the set of
coefficients $\mathcal{A}_{\beta,\alpha}$, so that \begin{equation}
    \mathcal{A}_{\beta,\alpha} = \{a_{ij}\} , i=
\alpha,\dots,\beta-2; j > i+1. \end{equation}

\begin{theorem}[]{The GLV equations are ($\beta, \alpha$)-reducible via IS if $a = 0\,\,
    \forall\,\, a \in \mathcal{A}_{\beta,\alpha}$, 
    and $a \neq 0\,\,\forall\,\, a \in \mathcal{A}_{\beta,\alpha}^c$.}
\end{theorem}
\begin{proof}
    In Lemma~\ref{lem:memSS-1}, we completed the reduction from $\beta$ to
    $\beta-1$ equations. The next step reduces to $\beta-2$. Let's start with the $\beta-1$ equations:
\begin{subequations}
    \begin{align}
        \dot x_1 &= b_1 x_1 + (a_{11} x_1  + \dots + a_{1\beta}\bm{\chi_{\beta}^{1:\beta-1}}) x_1 \label{eq:2S-1,1}\\
        \dot x_2 &= b_2 x_2 + (a_{21} x_1  + \dots + a_{2\beta}\bm{\chi_{\beta}^{1:\beta-1}}) x_2 \label{eq:2S-1,2}\\
        &\vdots \nonumber \\
        \dot x_{\beta-2} &= b_{\beta-2} x_{\beta-2} + \Bigl( a_{\beta-2,1} x_1  + \dots +
        a_{\beta-2,\beta-2}x_{\beta-2}+ a_{\beta-2,\beta-1}x_{\beta-1}+ a_{\beta-2,\beta}\bm{\chi_{\beta}^{1:\beta-1}}\Bigr)
        x_{\beta-2} \label{eq:2S-1,S-2}\\
        \dot x_{\beta-1} &= b_{\beta-1} x_{\beta-1} + \Bigl(a_{\beta-1,1} x_1 + \dots +
        a_{\beta-1,\beta}\bm{\chi_{\beta}^{1:\beta-1}}\Bigr)
        x_{\beta-1} \label{eq:2S-1,S-1}.
    \end{align}
\end{subequations}
Now we rewrite $x_{\beta-1}$ in terms of the remaining $\beta-2$ variables; again, we could use any of
the first $\beta-2$ equations to do so. Let's choose the ODE for $x_{\beta-2}$ with the assumption that
$a_{\beta-2,\beta} = 0$:
    \begin{align}
        x_{\beta-1} &= \frac{1}{a_{\beta-2,\beta-1}} \biggl(\frac{1}{x_{\beta-2}} (\dot x_{\beta-2} - b_{\beta-2}
        x_{\beta-2}) - a_{\beta-2,1} x_1 -\dots - a_{\beta-2,\beta-3}x_{\beta-3} \biggr) \equiv \bm{y_{\beta-1}^{1:\beta-2}}.
    \end{align}
Substituting $\bm{y_{\beta-1}^{1:\beta-2}}$ into \ref{eq:2S-1,S-1},
    \begin{align}
        x_{\beta-1} &= \int_0^t \left(b_{\beta-1} \bm{y_{\beta-1}^{1:\beta-2}} + (a_{\beta-1,1} x_1 + \dots +
        a_{\beta-1,\beta}\bm{\chi_{\beta}^{1:\beta-1}})  \bm{y_{\beta-1}^{1:\beta-2}}\right)  \bm{\chi_{\beta-1}^{1:\beta-2}} . \end{align}
    And finally, substituting $\bm{\chi_{\beta-1}^{1:\beta-2}}$ back into the first $\beta-2$ equations,
\begin{subequations}
    \begin{align}
        \dot x_1 &= b_1 x_1 + \left(a_{11} x_1 + a_{12}x_2 + \dots+a_{1\beta}\bm{\chi_\beta^{1:\beta-1}}\right)
        x_1 \label{eq:dx1a}\\
        \dot x_2 &= b_2 x_2 + \left(a_{21} x_1 + a_{22}x_2 +\dots+ a_{2\beta}\bm{\chi_\beta^{1:\beta-1}}\right)
        x_2 \label{eq:dx2a}\\
        &\vdots \nonumber\\
        \dot x_{\beta-2} &= b_{\beta-2} x_{\beta-2}  + \left(a_{\beta-2,1} x_1 + a_{\beta-2,2}x_2 +
        \dots+a_{\beta-2,\beta}\bm{\chi_{\beta-1}^{1:\beta-1}}\right) x_{\beta-2} \label{eq:dxS-2a}.
   \end{align}
\end{subequations}

This process can be repeated once more, yielding a system of $\beta-3$ equations, if we also assume
$a_{\beta-3,\beta}  = a_{\beta-3, \beta-1} = 0$. And we may find the equivalent system in $\beta-4$ equations if the set
of zero interaction terms also includes $a_{\beta-4,\beta}$, $a_{\beta-4,\beta-1}$, and $a_{\beta-4,\beta-2}$.
Continuing in this way, the model can be reduced from $\beta$ to $s$ equations if $a = 0 \,\, \forall a \in
    \mathcal{A}_{\beta,s}$. Note the assumption that $a \neq 0\,\,\forall\,\,a \in
    \mathcal{A}_{\beta, \alpha}^c$ allows for division by these coefficients. 
\end{proof}

Let us examine this set $\mathcal{A}_{\beta,\alpha}$ in more detail.  We compute the fraction $\rho$
of these ``zeroed'' terms out of all interaction coefficients.
\begin{corollary}
Let $\lambda = \alpha/\beta$. Then
$\lim_{\beta \to \infty} \rho = \frac{1}{2}(1 - \lambda)^2.$
\end{corollary}
%Visually, the matrix $\bm{A}$ must look like:
%\begin{equation} \begin{bmatrix} 
%    a_{11} & a_{12} & \hdots & a_{1S}\\
%    a_{11} & a_{12} & \hdots & a_{1S}\\
%    a_{11} & a_{12} & \hdots & a_{1S}\\
%    a_{11} & a_{12} & \hdots & a_{1S}
%\end{bmatrix}.\end{equation}
\begin{proof}
Given $\beta$ species in the original model, there are $\beta^2$ interaction terms, and
    $|\mathcal{A}_{\beta,\alpha}| = \frac{(k-1)k}{2}$. In terms of $\alpha$, $\beta$, then $\rho$
    is:
\begin{align}
    \rho &= \frac{(k-1)k}{2 \beta^2} \\
    &= \frac{(\beta - s - 1)(\beta-s)}{2\beta^2} \\
    &= \frac{\beta^2 + s^2 - 2\beta s -\beta -s}{2\beta^2}.
\end{align}
Then the fraction $\rho$ may be written as
\begin{align}\rho &= \frac{1}{2}(1 + \lambda^2 - 2\lambda) -
\frac{1}{2\beta}(1 + \lambda)\\
    &= \frac{1}{2}(1 - \lambda)^2 - \frac{1}{2\beta}(1 + \lambda)
.\end{align}
In the limit as $\beta \to \infty$ and for a fixed $\lambda$, we have
\begin{align}
\lim_{\beta \to \infty} \rho = \frac{1}{2}(1 - \lambda)^2.\end{align}
\end{proof}

\begin{figure}[ht]
 \centering
 \includegraphics[width=.5\textwidth]{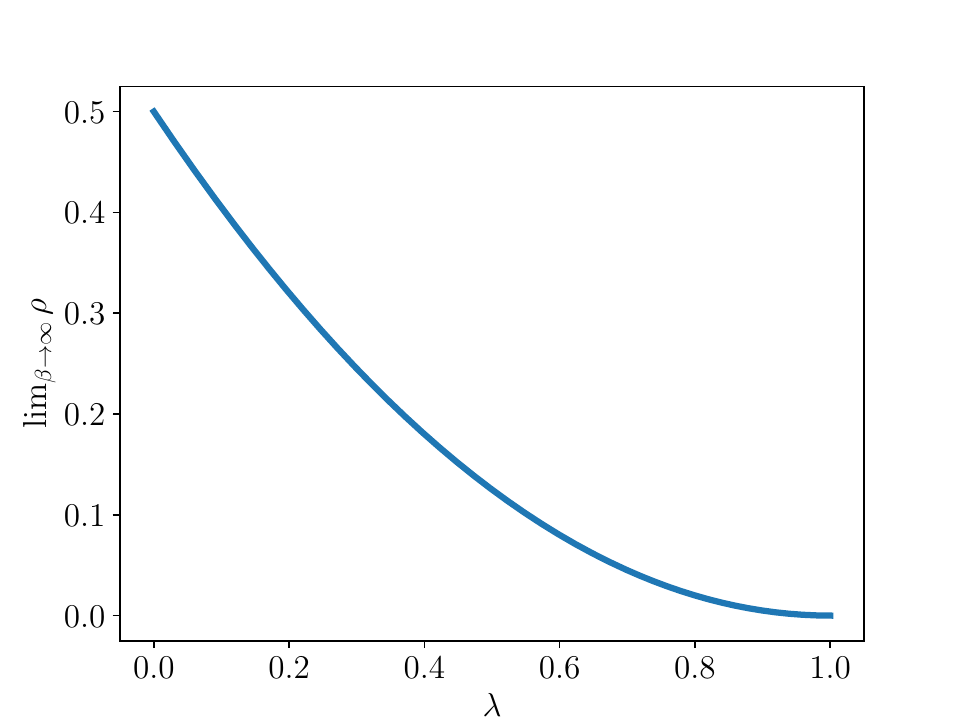}
 \caption{The fraction of interaction that must be zero for exact reduction decreases nonlinearly as
    $\alpha \to \beta$. \label{fig:lambda}}
\end{figure}
This limiting value is plotted for $\lambda \in [0,1]$ in figure~\ref{fig:lambda}. Note that in this
limit, when $\lambda=0$, \emph{i.e.,} there is complete reduction, then half of the interaction
terms vanish.  At the other extreme, when $\lambda = 1$, \emph{i.e.,} $\alpha=\beta$ so there is no
reduction, then of course all interaction terms can be nonzero.

Note that the  conditions on the set $\mathcal{A}_{\beta,\alpha}$ and its complement are not necessary, but
sufficient. They are not necessary at least in the sense of non-uniqueness: reordering the
variables or making different choices about which substitutions to make would lead to a different
set of zeroed interactions terms. In any case, we end this section with the following conjecture.

\begin{conj}[]For $(\beta,\alpha)$-reducibility, it is necessary that at least
$|\mathcal{A}_{\beta,\alpha}|$ terms are zero, where, for $\alpha = \beta-k$,
$|\mathcal{A}_{\beta,s}| = \frac{1}{2}k(k-1)$.  \end{conj}

%***************************************************************************************************
\section{Approximate model reduction}\label{sec:app} Motivated by the previous results, the final topic
of this paper explores approximate model reduction. In this section, we consider the situation that
information---interaction coefficients, growth rates, and some observations---is only available
about a subset of $\alpha$ out of $\beta$ species. We aim to create an approximate model for those
$\alpha$ species without including any information about the excluded $\beta-\alpha$ variables.

\subsection{Algebraic substitutions}\label{ssec:app-alg} In the algebraic method, eliminated variables are
exchanged for higher derivatives about the remaining variables. 

\begin{example}[Approximate AS for $\beta=2$, $\alpha=1$] The goal is to find an \textit{approximate} model
    for $x_1$, in terms of only $x_1$. 
\end{example}
The \textit{original} model is, as before:
\begin{subequations}
\begin{align}
    \dot x_1 &=  r_1  x_1 + ( a_{11} x_1 +  a_{12} x_2) x_1 \\%\label{eq:dx1}\\
    \dot x_2 &=  r_2  x_2 + ( a_{21} x_1 +  a_{22} x_2) x_2. %\label{eq:dx2}.
\end{align}
\end{subequations}

Motived by the AS above, we write $x_2$ as an expansion in $x_1$ and its higher derivatives:
\begin{equation} x_2 \rightarrow \kappa_0 x_1 + \kappa_1 \dot x_1 .\end{equation} 
    This leads to an approximate model for $x_1$:
    \begin{align} \dot x_1 &=r_1 x_1 + (a_{11}x_1 + a_{12}x_2)x_1  \nonumber \\
        &\approx r_1 x_1 + \left(a_{11}x_1 + a_{12}(\kappa_0 x_1 + \kappa_1 \dot x_1)\right)x_1 \nonumber \\ 
    &= r_1 x_1 + a_{11}x_1^2 + \delta_0 x_1^2 + \delta_1 \dot x_1 x_1, \end{align}
    where $\delta_0 = a_{11}\kappa_0$ and $\delta_1 = a_{11}\kappa_1$.

    It may also be of interest to see the behavior of $x_1$ given the remaining terms from the
    original model, i.e., $\dot x = r_1 x_1 + a_{11}x_1^2$. We call this the \textit{naive reduced} model.

Let us specify the coefficients of the original model as: \begin{equation} \bm{r} = \begin{bmatrix}5
    \\ 3\end{bmatrix}, \quad\quad  \bm{A} = \begin{bmatrix} -3 & -1 \\-1 & -2\end{bmatrix},
\end{equation} 
and set the initial conditions as $x_1(0) = 0.289$ and $x_2(0) = 0.665$. (These values are
approximate. The precise initial conditions were chosen randomly from a standard lognormal
distribution.)

\begin{figure}[ht]
 \centering
 \includegraphics[width=.5\textwidth]{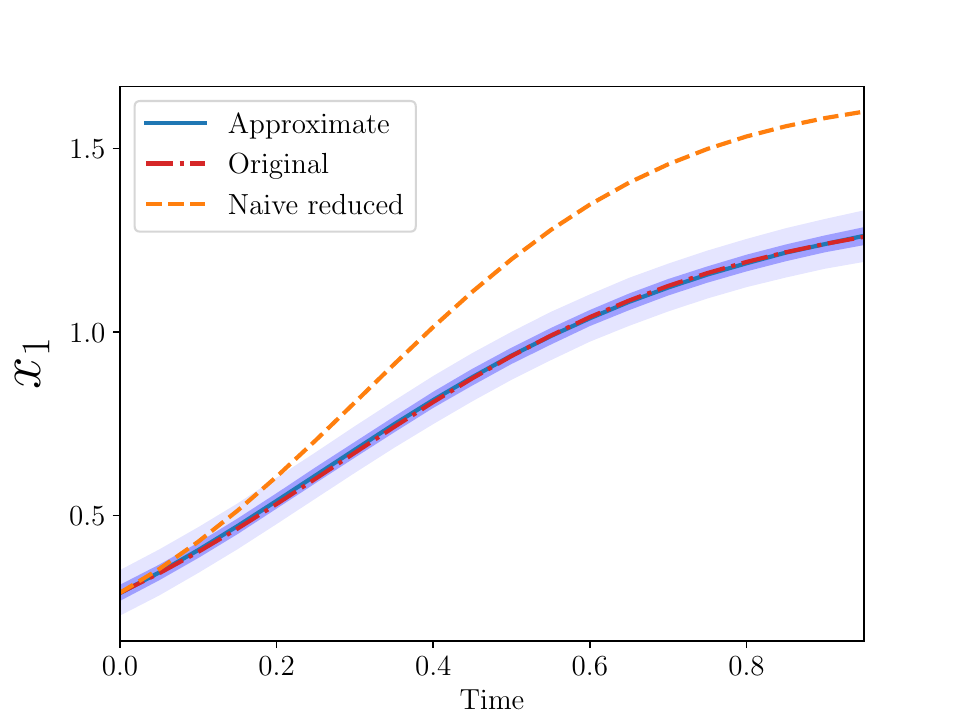}
 \caption{Approximate algebraic model for $\beta=2$, $\alpha=1$. The darker blue band represents the 50\%
    confidence interval (CI), and the lighter blue the 95\% CI.  \label{fig:2v1alg}}
\end{figure}

The introduced coefficients $\delta_0$ and $\delta_1$ are unknown and must be calibrated. We assume
data about $x_1$ over time. Calibration is performed under a Bayesian framework; the posterior means
are $\bar\delta_1 \approx -0.280 $ and $\bar \delta_2 \approx -0.357$.  Code to run Bayesian inverse
and forward problems is available
here: \texttt{github.libqueso.com} \cite{prudencio2012parallel}. %details are given in Appendix~\ref{app:cal}.  The
The trajectory of $x_1$ is shown in Figure~\ref{fig:2v1alg} from the three models: original, naive
reduced, and approximate.
The approximate model is quite close to the original. Importantly, the  approximate model
output completely covers the original model output, which lies entirely  within the approximate
model 50\% confidence interval.

\begin{example}[Approximate AS for $\beta=3$, $\alpha=1$.]\label{ex:as31} Here we try the same
    framework with a slightly larger reduction: from 3 to 1. In this case, note that we still
    replace $x_2$ and $x_3$ with our approximate expression, but still only two new coefficients
    need be introduced.
\end{example}

Recall the original model for $x_1$ is:
\[ \dot x_1 = r_1 x_1 + (a_{11} x_1 + a_{12}x_2 + a_{13}x_3) x_1.\]
The approximate model is then:
\begin{align} \dot x_1 &= r_1 x_1 + \Bigl(a_{11} x_1 + a_{12}(\kappa_1 x_1 + \kappa_2 \dot x_1) +
    a_{13}(\gamma_1 x_1 +\gamma_2 \dot x_1)\Bigr) x_1\nonumber\\
    &= r_1 x_1 + a_{11} x_1^2 + (a_{12}\kappa_1+ a_{13}\gamma_1) x_1^2 + (a_{12}\kappa_2 +
    a_{13}\gamma_2) x_1 \dot x_1\nonumber \\
&= r_1 x_1 + a_{11} x_1^2 + \delta_1 x_1^2 + \delta_2 x_1 \dot x_1\end{align}
where $\delta_1 = a_{12}\kappa_1+ a_{13}\gamma_1$ and $\delta_2 = a_{12}\kappa_2 +
    a_{13}\gamma_2$.
Thus, we can again calibrate simply $\delta_1$ and $\delta_2$.

Note the naive reduced model is again just
\[ \dot x_1 = r_1 x_1 + a_{11} x_1^2. \]

Specifically, we set
\begin{equation} \bm{r} = \begin{bmatrix}2.2 \\ 1.7 \\ 3.1\end{bmatrix}, \quad\quad  \bm{A} =
    \begin{bmatrix} -1.2 & -0.5 & -0.1 \\0.2 & -1.5 & -0.8 \\ -0.3 & -0.7 & -1.1\end{bmatrix},
\end{equation} 
and the initial conditions to $x_1(0) = 0.289$, $x_2(0) = 0.665$, and $x_3(0) = 5.01$.

\begin{figure}[ht]
 \centering
 \includegraphics[width=.5\textwidth]{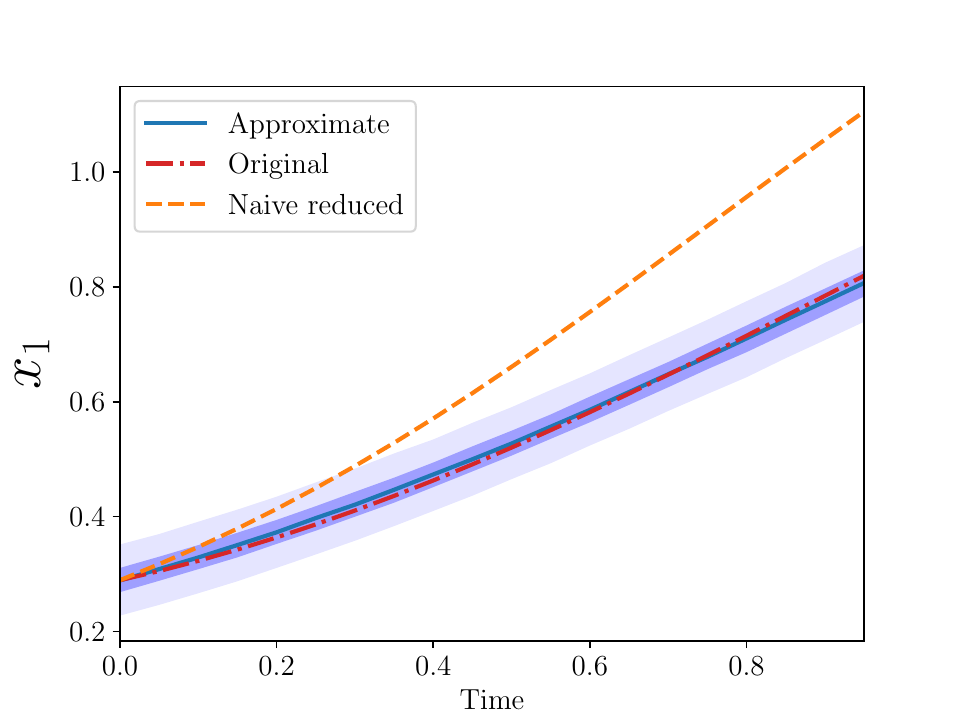}
 \caption{Approximate algebraic model for $\beta=3$, $\alpha=1$. \label{fig:3v1alg}}
\end{figure}

The posterior means are $\bar \delta_1 \approx -0.311$ and $\bar \delta_2 \approx -0.432$.
The trajectories from the three models are plotted in Figure~\ref{fig:3v1alg}.
Again, the approximate model is a good match to the original for $x_1$. Note that ten
coefficients from the original model are omitted (eight interaction and two growth rate
coefficients), while only two new parameters are introduced, $\delta_1$ and $\delta_2$.
%Investigating why such an exchange is even possible is beyond the scope of this paper, but will be
%studied in immediate future work.

\subsection{Integral substitutions}\label{ssec:app-mem}
Another approach is instead inspired by the IS above, so that information about
eliminated variables is replaced by memory, or integral, information of the remaining species.
\begin{example}[Approximate IS for $\beta=2$, $\alpha=1$]
Starting, as usual, with  the original formulation
    \[\dot x_1 =  r_1  x_1 + a_{11} x_1^2 +  a_{12} x_2 x_1\]
we aim for an approximate model by eliminating $x_2$.
\end{example}

In this example, we also try a slight variation from the algebraic one in that we replace the entire
term containing $x_2$:
\begin{equation}a_{12}x_1 x_2 \rightarrow \delta_1 x_1 + \delta_2 \int_0^t x_1(\tau) d\tau.\end{equation}
Moreover, since, in all of these examples, we assume some data about the species of interest, we can
approximate this integral directly from the data, which is generated according to the original model
(and not the $x_1$ as given by the approximate model). In these numerical results, this integral is
estimated numerically using a simple trapezoid rule. Call this numerical approximation
$\mathcal{I}(t)$:
\begin{equation} \mathcal{I}(t) \approx \int_0^t x^*_1(\tau) d\tau, \end{equation}
    where the $*$ in $x^*_1$ indicates that this integral is estimated from the data generated by
    the original model.
Then, we pose the approximate model for $x_1$ as
\begin{equation} \dot x_1 =   r_1  x_1 + a_{11} x_1^2 + \delta_1 x_1 + \delta_2 \mathcal{I}(t)
.\end{equation}

\begin{figure}[ht]
 \centering
 \includegraphics[width=.5\textwidth]{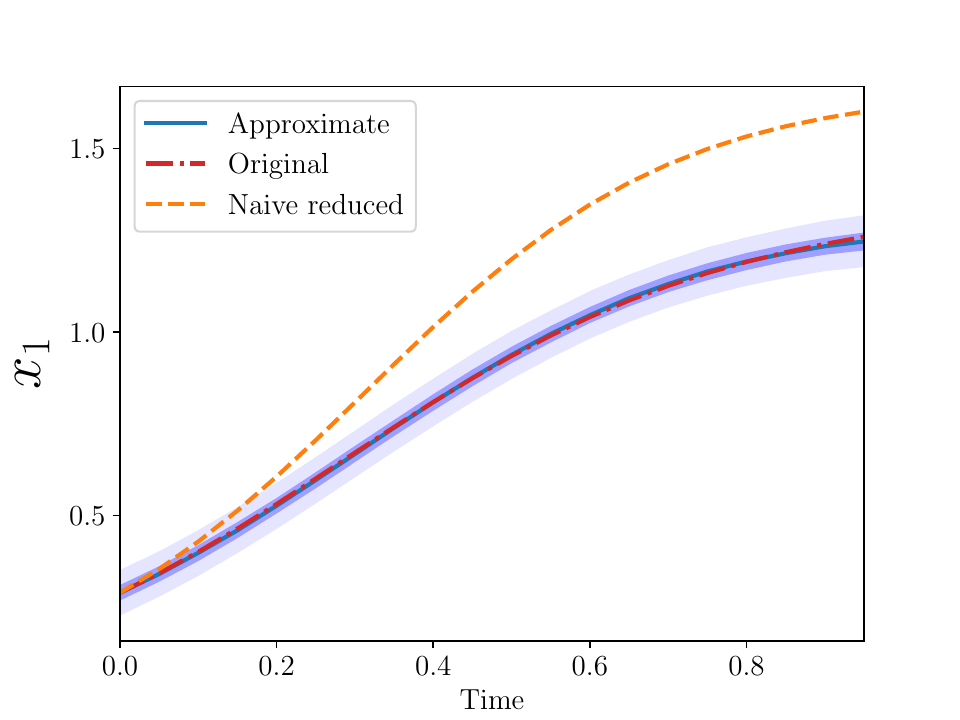}
 \caption{Approximate integral model for $\beta=2$, $\alpha=1$. \label{fig:2v1mem}}
\end{figure}

The posterior means in this case are $\bar \delta_1 \approx -0.778 $ and $\bar \delta_2 \approx -0.449 $.
The different trajectories for this example are shown in Figure~\ref{fig:2v1mem}.

\begin{example}[Approximate IS for $\beta=3$, $\alpha=1$]
    Here we present the final example, with the same original model as Example~\ref{ex:as31}
(approximate AS for $\beta=3$, $\alpha=1$), and the same approximate model as the previous example.
\end{example}
Here we find $\bar \delta_1 \approx -0.452$ and $\bar \delta_2 \approx -0.168$.

\begin{figure}[ht]
 \centering
 \includegraphics[width=.5\textwidth]{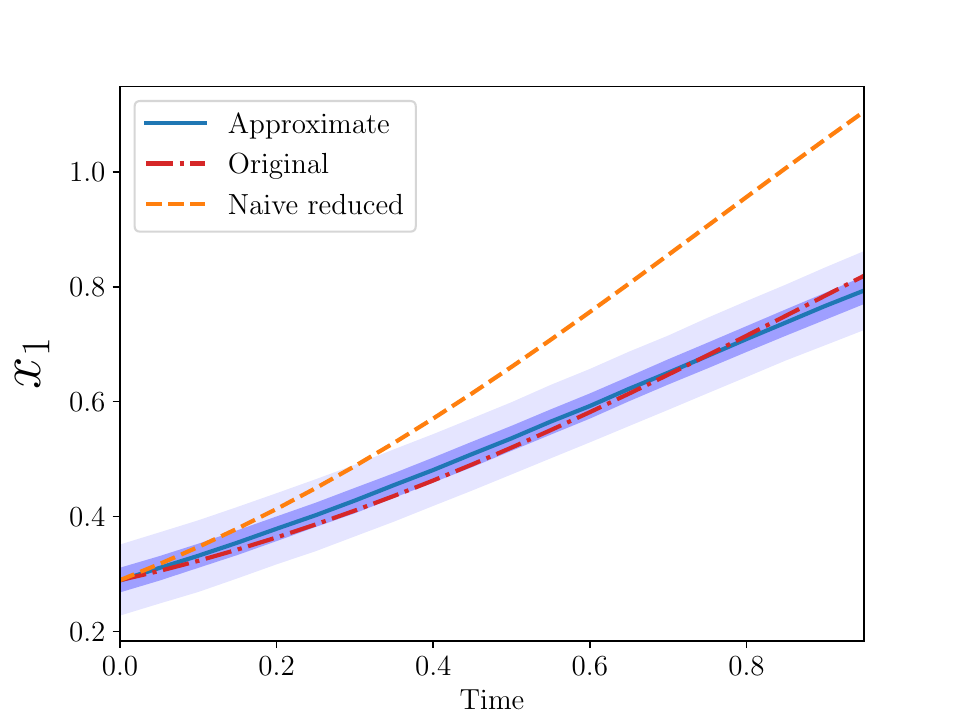}
 \caption{Approximate integral model for $\beta=3$, $\alpha=1$. \label{fig:3v1mem}}
\end{figure}

The results are shown in Figure~\ref{fig:3v1mem}.
The approximate integral models also show quite good agreement with the original model.

Note that each example here has preserved species correspondence: the species of interest is $x_1$,
and that is precisely what is modeled. This would readily generalize to approximate models for
$\alpha>1$.  And while the resulting differential equations now include either a derivative or an
approximate integral, the numerical solution is rather easily obtained without introducing much
numerical complexity. \cite{morrison2020data} provides an expanded investigation into 
similar AS-type approximate models.

%Note that when $A$ is symmetric, then $\rho = ( 1 - \alpha)^2$.
%***************************************************************************************************
\section{Discussion} \label{sec:con}
In summary, we can reduce a GLV system exactly by introducing integral terms of the remaining
species, or by introducing their higher derivatives. In a sense, some variables may be eliminated in
exchange for this extra information about the reduced set. A few interesting results emerge from
this work: First, a reduction from $\beta$ to $\beta-1$ variables is always possible. Second, a reduction
from $\beta$ to $\alpha$ is possible if the interaction terms in a specified set are zero. With
$\lambda = \alpha/\beta$, the fraction of terms in this set is approximately $(1-\lambda)^2/2$.
Third, we saw that the two different methods in the $\beta=2$ case yielded a single ODE for $x_1$,
but either by maintaining the structure of the original ODE for $x_1$, or for $x_2$.  Fourth, after
reducing to $\beta-1$ variables, an entanglement of species prevents further (exact) reduction
without breaking the original model in some other way.

In this paper, this break was achieved by setting interaction terms to zero. Of course, in some
systems, this may be a reasonable assumption. This work also serves as motivation or evidence of how
to approximate the role of the eliminated species using only the reduced set.  Preliminary results
suggest reduced models that are no longer exact, but still preserve species correspondence and are
computationally feasible.  Specific implementations and their effectiveness is an open and rich area
for future study.

\paragraph{Supplementary.} To run forward and inverse problems for the approximate models, code is
available here: \texttt{github.com/rebeccaem/enriched-glv} \cite{morrison2020glvcode}.

%%%%%%%%%%%%%%%%%%%%%%%%%%%%%%%%%%%%%%%%%%
\paragraph{Funding.} This research received no external funding.

%%%%%%%%%%%%%%%%%%%%%%%%%%%%%%%%%%%%%%%%%%
\paragraph{Acknowledgments.} I would like to acknowledge Tim Lyons, Sriram Sankaranarayanan, and Carl Simpson
for many helpful discussions about this work.

%%%%%%%%%%%%%%%%%%%%%%%%%%%%%%%%%%%%%%%%%%
%\conflictsofinterest{The author declares no conflict of interest.} 

%\externalbibliography{yes}
\bibliographystyle{plain}
\bibliography{citations.bbl}

\begin{thebibliography}{10}

\bibitem{barabas2016effect}
G.~Barab{\'a}s, M.~J.~Michalska-Smith, and S.~Allesina.
\newblock The effect of intra-and interspecific competition on coexistence in multispecies communities.
\newblock {\em The American Naturalist}, 188(1):E1--E12, 2016.

\bibitem{bayarri2007framework}
M.~J. Bayarri, J.~O. Berger, R.~Paulo, J.~Sacks, J.~A. Cafeo, J.~Cavendish, C.-H. Lin, and J.~Tu.
\newblock A framework for validation of computer models.
\newblock {\em Technometrics}, 49(2):138--154, 2007.

\bibitem{dantas2018calibration}
E.~Dantas, M.~Tosin, and A.~C. Jr.
\newblock Calibration of a {S}{E}{I}{R}–{S}{E}{I} epidemic model to describe the {Z}ika virus outbreak in {B}razil.
\newblock {\em Applied Mathematics and Computation}, 338:249 -- 259, 2018.

\bibitem{frassoldati2009simplified}
A.~Frassoldati, A.~Cuoci, T.~Faravelli, E.~Ranzi, C.~Candusso, and D.~Tolazzi.
\newblock Simplified kinetic schemes for oxy-fuel combustion.
\newblock In {\em 1st International conference on sustainable fossil fuels for future energy}, pages 6--10, 2009.

\bibitem{gear2003projective}
C.~W. Gear and I.~G. Kevrekidis.
\newblock Projective methods for stiff differential equations: Problems with gaps in their eigenvalue spectrum.
\newblock {\em Siam Journal on Scientific Computing}, 24(4):1091--1106, 2003.

\bibitem{givon2004extracting}
D.~Givon, R.~Kupferman, and A.~Stuart.
\newblock Extracting macroscopic dynamics: model problems and algorithms.
\newblock {\em Nonlinearity}, 17(6):R55--R127, aug 2004.

\bibitem{grilli2017feasibility}
J.~Grilli, M.~Adorisio, S.~Suweis, G.~Barabás, J.~R. Banavar, S.~Allesina, and A.~Maritan.
\newblock Feasibility and coexistence of large ecological communities.
\newblock {\em Nature Communications}, 8:8, 02 2017.
\newblock Copyright - Copyright Nature Publishing Group Feb 2017; Last updated - 2017-02-25.

\bibitem{haller2017exact}
G.~Haller and S.~Ponsioen.
\newblock Exact model reduction by a slow--fast decomposition of nonlinear mechanical systems.
\newblock {\em Nonlinear Dynamics}, 90(1):617--647, 2017.

\bibitem{harrington2017reduction}
H.~A. Harrington and R.~A. Van~Gorder.
\newblock Reduction of dimension for nonlinear dynamical systems.
\newblock {\em Nonlinear Dynamics}, 88(1):715--734, 2017.

\bibitem{jones1988global}
W.~Jones and R.~Lindstedt.
\newblock Global reaction schemes for hydrocarbon combustion.
\newblock {\em Combust. Flame; (United States)}, 73:3, 9 1988.

\bibitem{kazantzis2010new}
N.~Kazantzis, C.~Kravaris, and L.~Syrou.
\newblock A new model reduction method for nonlinear dynamical systems.
\newblock {\em Nonlinear Dynamics}, 59(1-2):183, 2010.

\bibitem{li1995global}
M.~Y. Li and J.~S. Muldowney.
\newblock Global stability for the {S}{E}{I}{R} model in epidemiology.
\newblock {\em Mathematical biosciences}, 125(2):155--164, 1995.

\bibitem{mezic2005spectral}
I.~Mezi{\'c}.
\newblock Spectral properties of dynamical systems, model reduction and decompositions.
\newblock {\em Nonlinear Dynamics}, 41(1-3):309--325, 2005.

\bibitem{morrison2020data}
R.~E. Morrison.
\newblock Data-driven corrections of partial lotka--volterra models.
\newblock {\em Entropy}, 22(11):1313, 2020.

\bibitem{morrison2020glvcode}
R.~E. Morrison.
\newblock rebeccaem/enriched-glv: Initial release, 2020.

\bibitem{morrison2018representing}
R.~E. Morrison, T.~A. Oliver, and R.~D. Moser.
\newblock Representing model inadequacy: A stochastic operator approach.
\newblock {\em SIAM/ASA Journal on Uncertainty Quantification}, 6(2):457--496, 2018.

\bibitem{murdoch1994continuum}
A.~I. Murdoch and D.~Bedeaux.
\newblock Continuum equations of balance via weighted averages of microscopic quantities.
\newblock {\em Proceedings of the Royal Society of London. Series A: Mathematical and Physical Sciences}, 445, 1994.

\bibitem{oberkampf2010verification}
W.~L. Oberkampf and C.~J. Roy.
\newblock {\em Verification and validation in scientific computing}.
\newblock Cambridge University Press, 2010.

\bibitem{oliver2015validating}
T.~A. Oliver, G.~Terejanu, C.~S. Simmons, and R.~D. Moser.
\newblock Validating predictions of unobserved quantities.
\newblock {\em Computer Methods in Applied Mechanics and Engineering}, 283:1310--1335, 2015.

\bibitem{pavliotis2008multiscale}
G.~Pavliotis and A.~Stuart.
\newblock {\em Multiscale Methods: Averaging and Homogenization}.
\newblock Springer-Verlag New York, 2008.

\bibitem{prudencio2012parallel}
E.~E. Prudencio and K.~W. Schulz.
\newblock {The parallel C++ statistical library ‘QUESO’: Quantification of Uncertainty for Estimation, Simulation and Optimization}.
\newblock In {\em Euro-Par 2011: Parallel Processing Workshops}, pages 398--407. Springer, 2012.

\bibitem{sahoo2017lie}
S.~Sahoo, G.~Garai, and S.~S. Ray.
\newblock Lie symmetry analysis for similarity reduction and exact solutions of modified {K}d{V}--{Z}akharov--{K}uznetsov equation.
\newblock {\em Nonlinear Dynamics}, 87(3):1995--2000, 2017.

\bibitem{yasuhiro1996global}
Y.~Takeuchi.
\newblock {\em Global dynamical properties of {L}otka-{V}olterra systems}.
\newblock World Scientific, 1996.

\bibitem{tartakovsky2011dimension}
A.~M. Tartakovsky, A.~Panchenko, and K.~F. Ferris.
\newblock Dimension reduction method for {O}{D}{E} fluid models.
\newblock {\em Journal of Computational Physics}, 230(23):8554 -- 8572, 2011.

\bibitem{wangersky1978lotka}
P.~J. Wangersky.
\newblock Lotka-{V}olterra population models.
\newblock {\em Annual Review of Ecology and Systematics}, 9:189--218, 1978.

\bibitem{williams2008detailed}
F.~A. Williams.
\newblock Detailed and reduced chemistry for hydrogen autoignition.
\newblock {\em Journal of Loss Prevention in the Process Industries}, 21(2):131 -- 135, 2008.
\newblock Hydrogen Safety.

\bibitem{zhdanov2002higher}
R.~Zhdanov.
\newblock Higher conditional symmetry and reduction of initial value problems.
\newblock {\em Nonlinear Dynamics}, 28(1):17--27, 2002.

\end{thebibliography}

\end{document}